\documentclass{amsart}
\usepackage{amsfonts,epsfig}

\newtheorem{theorem}{Theorem}
\newtheorem{proposition}[theorem]{Proposition}
\newtheorem{lemma}[theorem]{Lemma}
\newtheorem{corollary}[theorem]{Corollary}

 \newcounter{thlistctr}
 \newenvironment{thlist}{\
 \begin{list}%
 {\alph{thlistctr}}%
 {\setlength{\labelwidth}{2ex}%
 \setlength{\labelsep}{1ex}%
 \setlength{\leftmargin}{6ex}%
 \usecounter{thlistctr}}}%
 {\end{list}}

\newcommand{\RT}{{\sf RT}}

\newcommand{\bA}{{\bf A}}
\newcommand{\bB}{{\bf B}}

\newcommand{\bF}{{\bf F}}
\newcommand{\bG}{{\bf G}}

\newcommand{\bP}{{\bf P}}

\newcommand{\cK}{{\mathcal K}}
\newcommand{\cL}{{\mathcal L}}

\newcommand{\cP}{{\mathcal P}}

\newcommand{\ro}{{\rm o}}
\newcommand{\rO}{{\rm O}}

\newcount\minute        % Current minute within the hour
\newcount\hour          % Current hour (24-hour type)
\newcount\hourMins  % Temporary for taking hour modulo 60
\def\now%
{% Displays today's time as ``hh:mm''
  \minute=\time    % Number of minutes since midnight
  \hour=\time \divide \hour by 60 % Get hours
  \hourMins=\hour \multiply\hourMins by 60
  \advance\minute by -\hourMins % Hours modulo 60
  \zeroPadTwo{\the\hour}:\zeroPadTwo{\the\minute}%
}% --- \now ---
\def\timestamp%
{% Displays ``yyyy-mm-dd hh:mm''
  \today\ \now
}% --- \timestamp ---
\def\zeroPadTwo#1%
{% Left zero pad of the argument to 2 digits.  The argument
  \ifnum #1<10 0\fi    % Conditionally output a zero
  #1%    Then output the argument
}% --- \zeroPadTwo ---

\title {Sufficient Conditions for Labelled 0--1 Laws}

\author{Stanley Burris}
\address{Department of Pure Mathematics, University of Waterloo,
  Waterloo, Ontario, N2L 3G1, Canada}
\email{snburris@math.uwaterloo.ca}
\urladdr{www.thoralf.uwaterloo.ca}
\author{Karen Yeats}
\address{Department of Mathematics and Statistics, Boston University,
 111 Cummington St., Boston, MA 02215, USA}
\email{kayeats@bu.edu}


\begin{document}

\begin{abstract}
If $\bF(x) = e^{\bG(x)}$, where 
$\bF(x) = \sum f(n)x^n$ and $\bG(x) = \sum g(n)x^n$, with
$0\le g(n) =\rO\big(n^{\theta n}/n!\big)$,
$\theta \in (0,1)$, and 
$\gcd\big(n : g(n) >0\big)=1$, 
then $f(n)=\ro(f(n-1))$.

This gives an answer to Compton's request in Question 8.3 \cite{compton 1987} for 
an ``easily verifiable sufficient condition'' to show that an adequate class 
of structures has a labelled first-order 0--1 law, namely it suffices to 
show that the labelled 
component count function is $\rO\big(n^{\theta n}\big)$ for some
$\theta \in (0,1)$. It also provides the means to recursively construct an adequate class of structures with a labelled 0--1 law but not an unlabelled 0--1 law, answering Compton's Question 8.4.
\end{abstract}

\maketitle

\section{Introduction}
Exponentiating a power series can have the effect of smoothing out the behavior 
of the coefficients. In this paper we look at conditions
on the growth of the coefficients of $\bG(x) = \sum g(n)x^n$, where $g(n)\ge 0$,
which ensure that $f(n-1)/f(n) \rightarrow \infty$, where $\bF(x) = e^{\bG(x)}$.

Useful notation will be $f(n) \prec g(n)$ for $f(n)$ eventually less than $g(n)$ and $f(n) \in \RT_\infty$ for $f(n-1)/f(n) \rightarrow \infty$; the notation $\RT$ stands for the ratio test.

\section{The Coefficients of $e^{poly}$}\label{poly}

\begin{proposition}\label{hayman main}
Given
\begin{eqnarray*}
\bG(x)&:=&g(1)x + \cdots + g(d)x^d,\quad g(i) \ge 0,\ g(d) > 0,\\
&&\text{with }\gcd\big(j\le d : g(j)>0\big)\ =\ 1\\
\bF(x)&:=&\sum_{n\ge 0}f(n)x^n \ =\ e^{\bG(x)},
\end{eqnarray*}
the function $\bF(x)$ is Hayman-admissible. Thus 
\begin{equation} \label{a}
f(n)\ \sim\ \frac{\bF(r_n)}{{r_n}^n\cdot \sqrt{2\pi \bB(r_n)}}
\end{equation}
where $r_n$ is the unique positive solution to
\[
x\cdot \bG'(x) \ =\ n,
\]
and $ \bB(x) := x^2\bG''(x) + x\bG'(x)$.
\end{proposition}
\begin{proof}
Theorem X of Hayman \cite{hayman} shows that $\bF(x)$ is Hayman-admissible. Then the
rest of the claim is an immediate consequence of Corollary II 
of \cite{hayman} where
the saddle-point method is applied to find the asymptotics of the coefficients of
an admissible function. 
\end{proof}

\begin{corollary} \label{fin asymp}
For $\bF(x),\bG(x)$ as in the above proposition,
\begin{thlist}
\item
$f(n) \in\RT_\infty$,
\item
$\displaystyle f(n)\ =\ 
\exp\Big(-\frac{n\log n}{d}\big(1+\ro(1)\big)\Big).
$
\end{thlist}
\end{corollary}

\begin{proof}
Item (a) follows immediately from Corollary IV of Hayman \cite{hayman}.

For item (b) one uses $r_n\bG'(r_n)\ =\ n$ to obtain:
\begin{eqnarray*}
\Big(\frac{n}{cdg(d)}\Big)^{1/d}& \preceq& r_n\ \le
\Big(\frac{n}{dg(d)}\Big)^{1/d}\quad\text{for }c>1\\
r_n&=& 
\big(1+\ro(1)\big)
\Big(\frac{n}{dg(d)}\Big)^{1/d}\\
{r_n}^n& =& 
\big(1+\ro(1)\big)^n\Big(\frac{n}{dg(d)}\Big)^{n/d}\\
\bB(r_n)&=&\big(1+\ro(1)\big) d^2 g(d) \Big(\frac{n}{dg(d)}\Big)\ =\ 
\big(1+\ro(1)\big) dn\\
\bG(r_n)&=&\big(1+\ro(1)\big) g(d) {r_n}^d\ =\ 
\big(1+\ro(1)\big) \frac{n}{d}\\
\bF(r_n)&=&\exp\Big(\frac{n}{d} \big(1+\ro(1)\big)\Big).
\end{eqnarray*}
Apply these results to \eqref{a}.

\end{proof}

\section{Some Technical Lemmas}
Now we drop the assumption that $\bG(x)$ is a polynomial, but keep the
requirement
\begin{equation}\label{gcd is 1}
\gcd\big(n : g(n)>0\big)\ =\ 1.
\end{equation}
This implies that $f(n)\succ 0$.

Choose a positive integer $L\ge 2$ sufficiently large so 
\begin{eqnarray}
n>L\ \Rightarrow\ [x^n]\exp\Big(g(1)x + \cdots + g(L)x^L\Big)\ >\ 0. \label{bignuff2}
\end{eqnarray}

Given $\ell > L$ with $g(\ell)>0$ let
\begin{eqnarray}
\bG_0(x)&:=&\sum_{n \ge 1} g_0(n)x^n\ :=\ \sum_{1\le n\le \ell} g(n)x^n\nonumber\\
\bF_0(x)&:=&\sum_{n\ge 0} f_0(n)x^n\ :=\ \exp(\bG_0(x))\label{A0} \label{f0}\nonumber\\
\bG_1(x)&:=&\sum_{n\ge 1} g_1(n)x^n \ :=\ \sum_{n\ge \ell+1} g(n)x^n\nonumber\\
\bF_1(x)&:=&\sum_{n\ge 0} f_1(n)x^n\ :=\ \exp(\bG_1(x)).\label{A1}
\end{eqnarray}

\begin{lemma}\label{setup lemma}
Suppose $r\ge -1$ is such that
\begin{eqnarray}
ng(n)& =& \rO\big(f_0(n+r)\big). \label{cond on g}
\end{eqnarray}
Then
\[
nf_1(n) \ =\ \rO\big(f(n+r)\big). \label{f1 small}
\]
\end{lemma}

\begin{proof}
In view of \eqref{bignuff2} and \eqref{cond on g} we can choose $C_r$ such that
\begin{equation}\label{Cr}
ng(n) \ \le\ C_r f_0(n+r)\quad\text{for }n+r\ge L+1.
\end{equation}
Differentiating \eqref{A1} gives
\begin{eqnarray*}
nf_1(n)
& =& 
\sum_{j=\ell+1}^n j g(j)\cdot f_1(n-j)\\
&\le&
C_r\sum_{j=\ell+1}^n f_0(j+r)\cdot f_1(n-j)
\quad\text{by } \eqref{Cr} \\
&\le&
C_r\sum_{j=0}^{n+r} f_0(j)\cdot f_1(n+r-j)\\
&=& C_r f(n+r),
\end{eqnarray*}
the last line following from $\bF(x) = \bF_0(x)\cdot \bF_1(x)$.
\end{proof}

\begin{lemma} \label{better than monotone}
Suppose for every integer $r \geq -1$
\begin{eqnarray*}
ng(n)& =& \rO\big(f_0(n+r)\big). 
\end{eqnarray*}
Then $f(n-1)/f(n)\rightarrow \infty$.
\end{lemma}

\begin{proof}
Since $f_0(n) \in\RT_\infty$ by Corollary \ref{fin asymp} 
there is a monotone decreasing function
$\varepsilon(n)$ such that 
for any sufficiently large $M$ we have
$\varepsilon(n) > f_0(n)/f_0(n-1)$ for $n \ge M$, 
and $\varepsilon(n) \rightarrow 0$ as $n \rightarrow \infty$.

Thus
\begin{eqnarray*}
f(n)
&=&\sum_{0\le j\le n} f_0(j) f_1(n-j)\\
&=&\sum_{0\le j\le M-1} f_0(j) f_1(n-j)\ +\ \sum_{M\le j \le n} f_0(j) f_1(n-j)\\
&\le&\ro\big(f(n-1)\big)\ +\ \varepsilon(M)\sum_{M \le j\le n}
f_0(j-1) f_1(n-j)\\
&&\text{by Lemma \ref{setup lemma} and the choice of $\varepsilon$} \\
&\le&\ro\big(f(n-1)\big)\ +\ \varepsilon(M)f(n-1).
\end{eqnarray*}
Thus 
\[
\limsup_{n\rightarrow \infty} \frac{f(n)}{f(n-1)}\ \le\ \varepsilon(M),
\]
and as $M$ can be arbitrarily large it follows that
\[
\lim_{n\rightarrow \infty} \frac{f(n)}{f(n-1)}\ =\ 0.
\]
\end{proof}

\section{Main Result}
We are now in a position to prove the main result, making use of
\[
n!\ =\ \exp\big(n \log n \cdot (1+\ro(1)) \big),
\]
which follows from Stirling's result.

\begin{theorem} \label{main}
  Suppose $\bF(x) = \exp(\bG(x))$ with $\bF(x) = \sum_{n \geq
    0}f(n)x^n$, $\bG(x) = \sum_{n \geq 1}g(n)x^n$, and 
$f(n), g(n) \ge 0.$ 
Suppose also that $\gcd \big( n : g(n) > 0\big) = 1$
  and that 
for some $\theta \in (0,1)$
  \[
    g(n) = \rO\big(n^{\theta n}\big/n!\big).
  \]
Then
  \[
    f(n) \in\RT_\infty .
  \]
\end{theorem}

\begin{proof}
  From Corollary \ref{fin asymp}, 
for any integer $r \geq -1$ 
and any $\theta \in (0,1)$, by choosing $\ell>L$ such that $1/\ell < 1-\theta$,
we have
  \begin{eqnarray*}
    f_0(n+r) 
    &=& 
       \exp\Big(-\frac{(n+r)\log (n+r)}{\ell}\big(1+\ro(1)\big)\Big) \\
    &=& 
\exp\Big(-\frac{n\log n}{\ell}\big(1+\ro(1)\big)\Big) \\
    &\succ&\frac{n^{\theta n}}{(n-1)!}\ .
  \end{eqnarray*}
  Thus $ng(n) = \rO\big(f_0(n+r)\big)$.  The Theorem
  then follows from Lemma \ref{better than monotone}.
\end{proof}

\section{Best Possible Result}
The main result is in a natural sense the best possible.

\begin{proposition}\label{best possible}
Suppose $t(n)\ge 0$ with 
$\gcd\big(n : t(n)>0\big) = 1$ 
is such that for any $\theta \in (0,1)$ 
\[
t(n)\  \neq\  \rO\big(n^{\theta n}/n!\big).
\]
Then there is a sequence $g(n)\ge 0$ with 
$\gcd\big(n : g(n)>0\big) = 1$ and $g(n) \le t(n)$
but $f(n) \notin \RT_\infty$, where one has
$\bF(x) = \exp(\bG(x))$.
\end{proposition}

\begin{proof}
For $\theta \in (0,1)$ let
\[
S(\theta) \ =\ \big\{n \ge 1 : t(n) > n^{\theta n}/n!\big\}\, .
\]
Then $S(\theta)$ is an infinite set.

Let $M$ be such that $\gcd\big(n \le M : t(n) > 0\big) = 1$, and let
\begin{eqnarray*}
g_1(n)&:=& \begin{cases}
t(n)&\text{if $n \le M$}\\
0&\text{if $n>M$}
\end{cases}\\
\bG_1(x)&:=&\sum g_1(n)x^n\\
d_1&:=& \deg(\bG_1(x))\\
\bF_1(x)&:=&e^{\bG_1(x)}.
\end{eqnarray*}
 
For $m\ge 2$ we give a recursive procedure to define polynomials $\bG_m(x)$; 
then letting
\begin{eqnarray*}
d_m&:=& \deg(\bG_m(x))\\
\bF_m(x)&:=&e^{\bG_m(x)},
\end{eqnarray*}
by Proposition \ref{hayman main} 
\[
f_m(n)\ =\ \exp\Big(-\frac{n\log n}{d_m}\big(1+\ro(1)\big)\Big).
\]

To define $\bG_{m+1}(x)$, having defined $\bG_m(x)$, let
\[
h_m(n)\ :=\ \frac{1}{n!} \cdot n^{(1 - 1/2d_m)n} .
\]
Then
\[
\frac{h_m(n)}{f_m(n-1)}\ \rightarrow\ \infty\quad\text{as }n\rightarrow \infty.
\]
Thus we can choose an integer $d_{m+1}>d_m$ such that 
\begin{eqnarray*}
d_{m+1}&\in&S\Big(1-\frac{1}{2d_m}\Big)\\
h_m(d_{m+1})& >&f_m(d_{m+1}-1).
\end{eqnarray*}
This ensures that  $h_m(d_{m+1}) \le t(d_{m+1})$. Let 
\[
\bG_{m+1}\ :=\ \bG_m(x) + h_m(d_{m+1})x^{d_{m+1}}.
\]
Then
\[
\frac{f_{m+1}(d_{m+1})}{f_{m+1}(d_{m+1}-1)} \ 
\ge\ \frac{h_m(d_{m+1})}{f_m(d_{m+1}-1)}\ >\ 1.
\]

   Now let $\bG(x)$ be the nonnegative power series defined by the 
sequence of polynomials $\bG_m(x)$;  and let $\bF(x) = e^{\bG(x)}$. Then
$g(n) \le t(n)$ but $f(n) \notin \RT_\infty$ as
\[
\frac{f(d_{m+1})}{f(d_{m+1}-1)} \ =\ \frac{f_{m+1}(d_{m+1})}{f_{m+1}(d_{m+1}-1)}\ > \ 1.
\]

\end{proof}

\section{Application to 0--1 laws}

A class $\cK$ of finite relational structures is {\em adequate} if
it is closed under  disjoint union and the extraction of components.
One can view the structures as being {\em unlabelled} with the component
count function $p_U(n)$ and the total count function
$a_U(n)$, both counting up to isomorphism.
The corresponding {\em ordinary} generating series are
\[
\bP_U(x)\ :=\ \sum_{n\ge 1} p_U(n)x^n,\qquad
\bA_U(x)\ :=\ \sum_{n\ge 0} a_U(n)x^n
\]
connected by the fundamental equation
\begin{eqnarray} \label{fund un}
\bA_U(x)\ =\ \prod_{j\ge1} \big(1-x^j\big)^{-p_U(j)}.
\end{eqnarray}

One can also view the structures as being {\em labelled} (in all possible ways)
with the count
functions $p_L(n)$ for the connected members of $\cK$, 
and $a_L(n)$ for all members of $\cK$.
The corresponding {\em exponential} generating series are
\[
\bP_L(x)\ :=\ \sum_{n\ge 1} p_L(n)x^n/n!,\qquad
\bA_L(x)\ :=\ \sum_{n\ge 0} a_L(n)x^n/n!
\]
connected by the fundamental equation
\begin{eqnarray}  \label{fund lab}
\bA_L(x)\ =\ e^{\bP_L(x)}.
\end{eqnarray}

All references to Compton in this section are to the two papers
\cite{compton 1987} and \cite{compton 1989}.

\subsection{Unlabelled 0--1 Laws for Adequate Classes}
Let $\cK$ be an adequate class with unlabelled count functions and
ordinary generating functions as described above.
Compton showed that if the radius of convergence
$\rho_U$ of $\bA_U(x)$ is positive then $\cK$ has an unlabelled 0--1 
law\footnote{Given a logic $\cL$, $\cK$ has an unlabelled $\cL$ 0--1 law 
means that for any $\cL$ sentence $\varphi$,
the probability that $\varphi$ holds in $\cK$ will be either 0 or 1. 
In \cite{compton 1987}
Compton worked with first-order logic, in \cite{compton 1989} with monadic second-order 
logic. In both papers he simply used the phrases ``unlabeled 0--1 law'' and 
``labeled 0--1 law''.}
iff $a_U(n)\in\RT_1$, that is,
\[
\frac{a_U(n-1)}{a_U(n)}\ \rightarrow\ 1\quad\text{as } n\rightarrow \infty.
\]

$\cK$ is finitely generated if $r=\sum p_U(n) < \infty$, that is,
there are only finitely many connected structures in $\cK$. 
In the finitely generated case
the asymptotics for the coefficients $a_U(n)$ 
have long been known to have the simple polynomial 
form\footnote{This result is usually known as Schur's Theorem \cite[3.15.2]{wilf}.  One can easily find the asymptotics 
\eqref{change prob}
using a partial fraction decomposition of the right side of \eqref{fund un}.
The labelled case with finitely many components is more difficult---we needed 
to invoke Hayman's treatise \cite{hayman} just to obtain
the asymptotics for $\log a_L(n)/n!$ (see Corollary \ref{fin asymp}). 
}
\begin{equation}\label{change prob}
a_U(n)\ \sim\ Cn^{r-1}
\end{equation}
provided $ \gcd\big(n : p_U(n)>0\big)\, =\, 1.  $
Item \eqref{change prob} leads to the fact that $a_U(n)\in\RT_1$, and 
hence to an unlabelled 0--1 law.
In addition to using this result, Compton notes that
the work of Bateman and Erd\"os \cite{bateman erdos} shows that if
$p_U(n)\in\{0,1\}$, for all $n$, then one has $a_U(n)\in\RT_1$. 

Both of these results were subsumed in the powerful result of Bell \cite{bell 2002} 
which
says that if $p_U(n)$ is polynomially bounded, that is, there is a $c$ such that
$p_U(n) = \rO(n^c)$, then $a_U(n) \in \RT_1$.

\subsection{Labelled 0--1 Laws}

Compton shows that if $\rho_L$, the radius of convergence of $\bA_L(x)$,
is positive, then $\cK$ has a labelled 0--1 law iff
\begin{equation}\label{lab test}
\frac{a_L(n-k)/(n-k)!}{a_L(n)/n!}\rightarrow \infty\quad\text{whenever } p_L(k)>0.
\end{equation}
In particular
it suffices to show that $a_L(n)/n! \in \RT_\infty$.

Compton's method to show that a given adequate class of finite relational 
structures $\cK$ has a labelled 0--1 law is to 
show that its exponential generating function 
$\bA_L(x) = \sum a_L(n)x^n/n!$ is Hayman-admissible with 
an infinite radius of convergence. This guarantees that 
$a_L(n)/n! \in \RT_\infty$ (\cite{hayman}, Corollary IV). However, as Compton notes, 
showing that $\bA_L(x)$ is Hayman-admissible can be quite a challenge.

Question 8.3 of \cite{compton 1987} first asks if, in the {\em unlabelled} case, 
the result of Bateman and Erd\"os, namely $p_U(n) \in \{0,1\}$ implies 
$a_U(n) \in \RT_1$, 
can be extended to the much more general statement that $p_U(n) = \rO(n^k)$
implies $a_U(n)\in\RT_1$, yielding an unlabelled 0--1 law. As mentioned
earlier, this was proved to be true by Bell. 
The second part of Question 8.3 asks if there is a simple sufficient condition 
along similar lines for the labelled case. We can now answer this in the 
affirmative with a result that is an excellent parallel to Bell's result for
unlabelled structures.

\begin{theorem} \label{01 laws}
If $\cK$ is an adequate class of structures with 
\[
p_L(n)\ =\ \rO\Big(n^{\theta n}\Big)\quad \text{for some }\theta \in (0,1)
\]
then $a_L(n)/n! \in \RT_\infty$, and consequently $\cK$ has a 
labelled monadic second-order 0--1 law.
\end{theorem}

\begin{proof}
This is an immediate consequence of Theorem \ref{main} and 
Compton's proof 
that $a_L(n)/n! \in \RT_\infty$ guarantees such a 0--1 law.
\end{proof}

Now we list the examples of classes $\cK$ which Compton 
shows have a labelled 0--1 law,
giving $p_L(n)$ in each case. It is trivial to check in each case that $p_L(n) = \rO\Big( n^{n/2}\Big)$; thus the 0--1 law in each case follows from our Theorem \ref{01 laws}.
\begin{thlist}
\item {\sf 7.1 Unary Predicates}\quad $p_L(n) = 0$ for $n>1$.
\item {\sf 7.12 Forests of Rooted Trees of Height 1}\quad $p_L(n) = n$. 
\item {\sf 7.15 Only Finitely Many Components}\quad $p_L(n)$ is eventually 0.
\item {\sf 7.16 Equivalence Relations}\quad $p_L(n) = 1$. 
\item {\sf 7.17 Partitions with a Selection Subset}\quad  $p_L(n) = 2^n-1$. 
\end{thlist}
We can now augment this list by, in each case, coloring the members of $\cK$ by
a fixed set of $r$ colors in all possible ways. This will increase the 
original $p_L(n)$ by a factor of at most $r^n$.  This will still give $p_L(n) = \rO\big(n^{n/2}\big)$. 
Furthermore, in each of these colored cases let $\cP$ be any subset of the 
connected members, and let $\cK$ be the closure of $\cP$ under disjoint union. 
Each such $\cK$ has a labelled 0--1 law.

Another application of Theorem \ref{01 laws} is to answer Question 4 of \cite{compton 1987}
by exhibiting an adequate class $\cK$ such that 
$p_L(n) =\rO\big(n^{3n/4}\big) $, 
hence there is a labelled 0--1 law for $\cK$;
but also such that $\rho_U \in (0,1)$, so $\cK$ does not have an unlabelled 0--1 law.

\begin{figure}\label{broom}
\epsfig{file=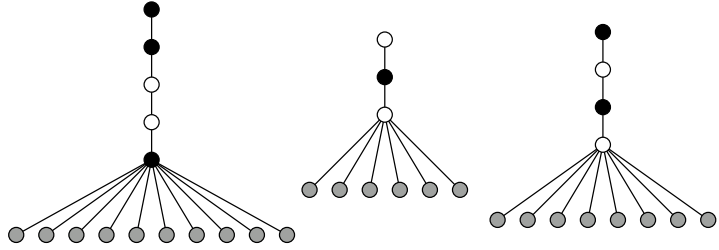, scale=1.0}
\caption{Brooms with two-colored handles}
\end{figure}

Let the components of $\cK$
 be the one-element tree $T_1$ along with rooted
trees $T_{3n}$ of size $3n$ and height $n$ consisting
of a chain $C_n$ of $n$ nodes, with an antichain $L_{2n}$ of $2n$ nodes 
(the leaves of the tree) below the least member of the chain; and the chain
$C_n$ is two-colored while the remaining nodes are uncolored. One can visualize these as brooms with 2-colored handles, see Figure \ref{broom}.

The number of unlabelled components is given by $p_U(1) = 1$, $p_U(3n) = 2^n$.
Thus the radius of convergence of the ordinary generating function
of $\cK$ is $\rho_U=\sqrt[3]{2}$. Since this
is positive and not 1 it follows from Theorem 5.9(ii) of \cite{compton 1987} that 
$\cK$ does not have an unlabelled 0--1 law.

For the number $p_L(3n)$ of labelled components of size $3n$: 
\begin{eqnarray*}
p_L(3n)&\le&2^n {3n \choose n} n!\\
&\le&2^n(3n)^n \exp\big(n\log n \cdot (1+\ro(1))\big)\\
&=&\exp\big(2n\log n \cdot (1+\ro(1))\big)\\
&=&(3n)^{(2/3)(3n)\big(1+\ro(1)\big)}\\
&=&\rO\Big((3n)^{(3/4)(3n)}\Big).
\end{eqnarray*}

Thus $p_L(n) = \rO\big(n^{3n/4}\big)$, so
$a_L(n)/n! \in \RT_\infty$ by Theorem \ref{01 laws}, showing that $\cK$ has a labelled 0--1 law.

\end{document}